\documentclass[11pt,a4paper]{article}

\usepackage[austrian,english]{babel}
\usepackage[latin1]{inputenc}
\usepackage[a4paper]{geometry}
\usepackage{makeidx}

\usepackage{hyperref}
\usepackage[pdftex]{graphicx}

\usepackage{amsmath}
\usepackage{amssymb}

\newenvironment{myproof}{\textnormal{\textbf{Proof.}}}{\hfill$\blacksquare$}
\newcommand{\defeq}{\mathrel{\mathop:}=}
\newcommand{\op}[1]{\operatorname{#1}}
\newtheorem{mycor}{Corollary}[section]
\newtheorem{mylem}{Lemma}[section]
\newtheorem{myprop}{Proposition}[section]
\newtheorem{mythrm}{Theorem}[section]

\begin{document}
\title{\textbf{Symmetry properties of the \\Novelli--Pak--Stoyanovskii algorithm}}
\author{Robin Sulzgruber
\thanks{
Supported by the Austrian Science Foundation FWF, grant S50-N15 in the framework of the Special Research Program ``Algorithmic and Enumerative Combinatorics'' (SFB F50).\newline
Email: \texttt{robin.sulzgruber@univie.ac.at}}
\\ Fakult\"at für Mathematik, Universit\"at Wien,\\ Oskar-Morgenstern-Platz 1, A1090 Vienna, Austria.}

\date{\today}
\maketitle
\begin{abstract}
The number of standard Young tableaux of a fixed shape is famously given by the hook-length formula due to Frame, Robinson and Thrall. A bijective proof of Novelli, Pak and Stoyanovskii relies on a sorting algorithm akin to jeu-de-taquin which transforms an arbitrary filling of a partition into a standard Young tableau by exchanging adjacent entries. Recently, Krattenthaler and M\"uller defined the complexity of this algorithm as the average number of performed exchanges, and Neumann and the author proved it fulfils some nice symmetry properties. In this paper we recall and extend the previous results and provide new bijective proofs.
\end{abstract}

\section{Introduction}

The main motivation for this work was a conjecture by Krattenthaler and M\"uller stating that the average runtime of the Novelli--Pak--Stoyanovskii algorithm remains the same whether it is applied row-wise or column-wise. Equivalently, the sorting of all fillings of a partition $\lambda$ requires the same number of exchanges as the sorting of all fillings of its conjugate $\lambda'$.

We consider two additional combinatorial objects, the drop function and exchange numbers, by means of which the complexity can be expressed. The aim of this paper is to provide formulae for the drop function (Theorem \ref{thrm:drop}), the exchange numbers (Theorem \ref{thrm:exchange}) and the complexity (Theorem \ref{thrm:complexity}). The symmetry in $\lambda$ and $\lambda'$ will follow easily from these results (Corollary \ref{cor:symm}). 

It should be mentioned that Neumann and the author obtain the symmetry result in a slightly more general setting, defining a family of closely related algorithms and proving that two such algorithms share the same complexity, drop function and exchange numbers if both algorithms produce each standard Young tableau equally often as an output. However, restricting ourselves to the original algorithm, we are able to deduce more explicit results. 

Finally, some of these new formulae generalise the hook-length formula, and demand bijective proofs. We succeed in providing bijections for the drop function (Theorem \ref{thrm:dropbij}) and the exchange numbers (Proposition \ref{prop:exchangebij}), and give a bijective version of the symmetry of the complexity at the end of Section \ref{sec:bijections}.

\bigskip
In the rest of this section we define the combinatorial objects we are going to consider, and recall the Novelli--Pak--Stoyanovskii algorithm.

\begin{figure}
\begin{center}
 \includegraphics[height=15mm]{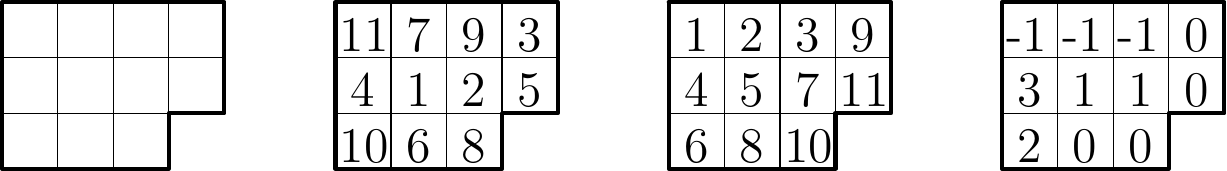}
 \caption{The partition $\lambda=(4,4,3)$, and a tabloid, a standard Young tableau and a hook tableau of shape $\lambda$.}
 \label{fig:partition}
\end{center}
\end{figure}

Let $n \in \mathbb N$. A partition of $n$ is a sequence $\lambda = (\lambda_1, \lambda_2, \dots)$ of non-negative integers such that $\lambda_i \geq \lambda_{i+1}$ and $\sum \lambda_i = n$. If this is the case we write $\lambda \vdash n$. For us the appropriate way to represent a partition is via its Young diagram $\lambda = \{ (i,j) \in \mathbb Z^2 : 1 \leq i, 1 \leq j \leq \lambda_j\}$. The elements $x = (i, j) \in \lambda$ are called the cells of the partition. We visualise a partition as a left-justified array with $\lambda_i$ cells in the $i$-th row. Thus, the cells are arranged like the entries of a matrix (see Figure \ref{fig:partition}). The conjugated partition of $\lambda$ is defined as
\begin{align*}
\lambda'
= \{(i, j) \in \mathbb Z^2 : 1 \leq i, 1 \leq j \leq \lambda_j'\}
\defeq \{(i, j) : (j, i) \in \lambda \}.
\end{align*}

Given a cell $x = (i, j)$ in $\lambda$ we define its arm $\operatorname{arm}_{\lambda}(x) \defeq \lambda_i - j$ as the number of cells to the right of $x$. The leg $\op{leg}_{\lambda}(x) \defeq \lambda'_j - i$ is defined as the number of cells below $x$. Furthermore, we define the hook-length of $x$ as $h_{\lambda}(x) \defeq \op{arm}_{\lambda}(x) + \op{leg}_{\lambda}(x) + 1$, and the cohook-length as $h'(x) \defeq i + j - 2$. Finally, we denote the set of top and left neighbours of $x$ by $N_{\lambda}^{-}(x) \defeq \{ (i-1,j), (i,j-1) \} \cap \lambda$ and the set of right and bottom neighbours of $x$ by $N_{\lambda}^{+}(x) \defeq \{ (i, j+1), (i+1, j)\} \cap \lambda$.


We want to consider integer fillings of partitions, i.e.\ maps $T: \lambda \to \mathbb Z$. The image $k = T(x)$ is called the entry of the cell $x$. The partition $\lambda$ is also called the shape of the filling. A tabloid is a bijection $T: \lambda \to \{1, \dots, n\}$, i.e.\ a filling where each entry between $1$ and $n$ occurs exactly once. We say a tabloid $T$ is sorted at $x$ if $T(x)$ is less than all entries among the bottom and right neighbours of $x$. A tabloid which is sorted at every cell is called a standard Young tableau. Finally, a hook tableau is a map $H: \lambda \to \mathbb Z$ such that $-\op{leg}_{\lambda}(x) \leq H(x) \leq \op{arm}_{\lambda}(x)$ for all $x \in \lambda$.

Let $\op{T}(\lambda)$ denote the set of all tabloids, $\op{SYT}(\lambda)$ the set of all standard Young tableaux and $\op{H}(\lambda)$ the set of all hook tableaux of shape $\lambda$. It is obvious that
\begin{align}\label{eq:TH}
\#\op{T}(\lambda) = n! \quad\text{ and }\quad
\#\op{H}(\lambda) = \prod_{x \in \lambda} h_{\lambda}(x).
\end{align}

The symmetric group $\mathfrak S_n$ acts simply transitively on $\op{T}(\lambda)$ via $\sigma(T) \defeq \sigma \circ T$. Let $T \in \op{T}(\lambda)$ and $\sigma = (k_0, \dots, k_r)$ be a cycle such that $T^{-1}(k_{i})$ is a right or bottom neighbour of $T^{-1}(k_{i-1})$ for all $1 \leq i \leq r$. Then $\sigma$ is called a forward slide on $T$ at $x$ if $T(x) = k_0$, $k_{i} < k_{i+1}$ for all $1 \leq i < r$ and $k_0 > k_r$. If instead $k_0 = T(x)$ and $k_{i} < k_{i+1}$ for all $0 \leq i < r$ then the inverse $\sigma^{-1}$ is called a backward slide on $T$ at $x$. For example $(7,1,2,5)$ is a forward slide on the (left) tabloid in Figure \ref{fig:partition} while $(1,2,8)^{-1}$ is a backward slide.

An exchange on $T$ is a transposition $\tau = (k,l)$ such that $k$ and $l$ are entries of neighbouring cells in $T$. If $\sigma_1 = (k_0, k_1, \dots, k_r)$ is a forward slide on $T$ and $\sigma_2 = (k_0, l_1, \dots, l_s)$ is a forward slide on $\sigma_1 \circ T$ then $\sigma_2 \sigma_1 = (k_0, k_1, \dots, k_r, l_1, \dots, l_s)$ is a forward slide on $T$ if and only if $k_r < l_1$. In particular $\sigma_2 \sigma_1$ has to be a forward slide on $T$ if $T$ is ordered at $T^{-1}(k_r)$. Conversely, every forward slide $\sigma = (k_0, k_1, \dots, k_{r+s})$ has a unique decomposition $\sigma = \sigma_2 \sigma_1$ into a forward slide $\sigma_1$ on $T$ of length $r$ and a forward slide $\sigma_2$ on $\sigma_1 \circ T$ of length $s$. We say $\sigma$ is an extension of $\sigma_1$, and call a forward slide maximal if it possesses no nontrivial extension. Along the same lines, each forward slide of length $r$ has a unique decomposition into $r$ exchanges
\begin{align*}
(k_0,\dots,k_r) = (k_r,k_0) \cdots (k_2, k_0)(k_1, k_0),
\end{align*}
where each transposition $(k_i, k_0)$ is an exchange on $(k_{i-1}, k_0) \cdots (k_1, k_0) \circ T$.


\bigskip
The number of standard Young tableaux is famously given by the hook-length formula due to Frame, Robinson and Thrall \cite{FRT1954}:
\begin{align} \label{eq:HLF}
f_{\lambda}
\defeq \#\op{SYT}(\lambda)
= \frac{n!}{\prod_{x \in \lambda} h_{\lambda}(x)}
\end{align}
While there are many proofs of this formula, we are particularly interested in the proof due to Novelli, Pak and Stoyanovskii \cite{NPS1997}. They construct a bijection
\begin{align} \label{eq:nps}
\Phi: \op{T}(\lambda) \to \op{H}(\lambda) \times \op{SYT}(\lambda),
\end{align}
hence, the respective cardinalities are equal and the formula follows from \eqref{eq:TH}. The induced map $\op{T}(\lambda) \to \op{SYT}(\lambda)$ which transforms each tabloid $T$ into a standard Young tableau, is especially nice since it is given as a sequence of maximal forward slides and can be fully explained without considering hook tableaux at all:

Impose the reverse lexicographic order on the cells of $\lambda$, i.e.\ $(i,j) \prec (k,l)$ if either $j < l$ or $j = l$ and $i < k$. Initialise $y$ to be the maximal cell with respect to $\prec$ and $T_y \defeq T$. Then iterate the following steps:

Given a pair $(y, T_y)$ let $x$ be precursor of $y$ with respect to $\prec$. Define a maximal forward slide $\sigma_x$ on $T_y$ at $x$ by letting $k_0 \defeq T(x)$ and choosing each $k_i$ minimal for $i \geq 1$. Set $T_x \defeq \sigma_x \circ T_y$ and return the pair $(x, T_x)$.

\begin{figure}
 \begin{center}
  \includegraphics[height=15mm]{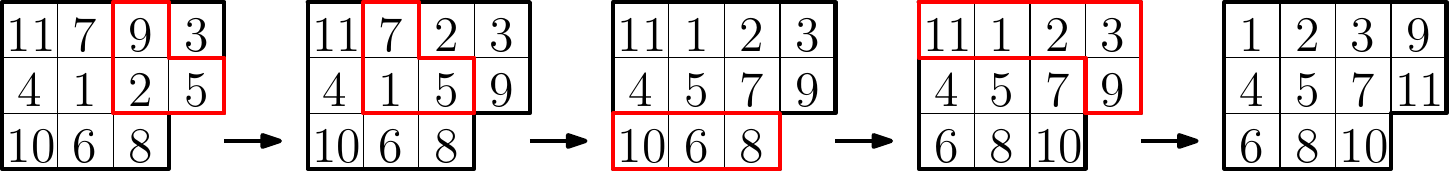}
  \caption{The non-trivial forward slides are $(9,2,5)$, $(7,1,5)$, $(10,6,8)$ and $(11,1,2,3,9)$.}
  \label{fig:sort}
 \end{center}
\end{figure}

The algorithm terminates after the minimal cell $(1,1)$ is reached. At each step the tabloid $T_x$ will be ordered at all cells $y$ with $y \succ x$. Thus, $T_{(1,1)}$ is a standard Young tableau. An example is given in Figure~\ref{fig:sort}.

The inverse algorithm is given as a sequence of backward slides, however, to determine the correct backward slide at each step, the hook tableau has to be taken into account. For details see \cite{NPS1997,Sagan2001}.

We remark that Neumann and the author define analogous sorting algorithms in \cite{NeuSul2013}, by replacing $\prec$ with a more general order $\prec_U$. This order is induced by an arbitrary standard Young tableaux $U$ of the same shape as $T$, and defined by $x \prec_U y$ if and only if $U(x) < U(y)$.

The order $\prec$ is special in the following sense: If we sort all tabloids in $\op{T}(\lambda)$ with respect to $\prec$ then each standard Young tableau in $\op{SYT}(\lambda)$ is met exactly $\prod_{x\in\lambda} h_{\lambda}(x)$ times. If we sort with respect to $\prec_U$ for an arbitrary standard Young tableau $U$, this need not be the case anymore. Interestingly, it turns out that if two algorithms arising from different orders, produce each standard Young tableau equally often as an output, then they must automatically share a lot of properties \cite[Corollary 4.5, Corollary 5.4]{NeuSul2013}. In this paper we restrict ourselves to the order $\prec$ and improve the previous general results for this special case.


\bigskip
The following property of the Novelli--Pak--Stoyanovskii algorithm is so useful that it deserves its own lemma.

\begin{mylem}[Invariance]\label{lem:invariance} Let $k, n \in \mathbb N$ such that $k \leq n$, let $\lambda \vdash n$ be a partition, $\pi \in \mathfrak S_n$ be a permutation that fixes all $l$ where $1 \leq l < k$, and $T \in \op{T}(\lambda)$ be a tabloid. Set $\tilde T \defeq \pi \circ T$ and let $T_x, \tilde T_x$ be the tabloids which arise during the sorting of $T$ and $\tilde T$, respectively. Then each permutation $\pi_x$ defined by $\tilde T_x = \pi_x \circ T_x$ fixes all $l$ with $1 \leq l < k$.
\end{mylem}

\begin{myproof} We apply induction on $x$ with respect to $\prec$. The claim is trivial for the maximal cell. Thus, let $x \in \lambda$ be a cell with successor $y \succ x$, and suppose that $\pi_y$ fixes all $l$ with $l < k$.

Let $\sigma_x = (a_0, a_1, \dots a_r)$ and $\tilde \sigma_x = (b_0, b_1, \dots, b_s)$ be the maximal forward slides defined as above, such that $a_0 = T_y(x)$ and $b_0  = \tilde T_y(x)$. Since we always prefer smaller entries over larger ones, and since all entries $l$ with $l < k$ appear in the same positions in $T_y$ and $\tilde T_y$, whenever we choose an entry $a_i < k$ we must make the same choice for $b_i$. Thus, we may decompose $\sigma_x$ and $\tilde \sigma_x$ into forward slides
\begin{align*}
\sigma_x = \pi_1 \circ (a_0, a_1, \dots, a_t) \qquad\text{ and }\qquad \tilde \sigma_x = \pi_2 \circ (b_0, a_1, \dots, a_{t}),
\end{align*}
such that $a_t < k$ and where $\pi_1$ and $\pi_2$ fix all $l$ with $l < k$. If $a_0 < k$ then $a_0 = b_0$ and $\sigma_x = \tilde \sigma_x$ fixes all $l$ with $l \geq k$. Hence,
\begin{align*}
\tilde T_x
= \sigma_x \circ \tilde T_y
= \sigma_x \pi_y \circ T_y
= \pi_y \sigma_x \circ T_y
= \pi_y \circ T_x.
\end{align*}
On the other hand, if $a_0 \geq k$ then also $b_0 \geq k$. Write $\pi_y = \pi_3 \circ (a_0, b_0, c_1, \dots, c_p)$ such that $\pi_3$ fixes $a_0$, $b_0$ and all $l$ with $l < k$. Then we compute
\begin{align*}
\tilde T_x
&= \pi_2 (b_0, a_1, \dots, a_t) \circ \tilde T_y \\
&= \pi_2 (b_0, a_1, \dots, a_t) \pi_3 (a_0, b_0, c_1, \dots, c_p) \circ T_y \\
&= \pi_2 \pi_3 (b_0, a_1, \dots, a_t) (a_0, b_0, c_1, \dots, c_p) \circ T_y \\
&= \pi_2 \pi_3 (a_0, b_0, c_1, \dots, c_p) (a_0, a_1, \dots, a_t) \circ T_y \\
&= \pi_2 \pi_y \pi_1^{-1} \pi_1 (a_0, a_1, \dots, a_t) \circ T_y
= \pi_2 \pi_y \pi_1^{-1} \circ T_x.
\end{align*}
Thus, we have $\pi_x = \pi_y$ in the first case and $\pi_x = \pi_2  \pi_y \pi_1^{-1}$ in the second case.
\end{myproof}

\bigskip
Krattenthaler and M\"uller were interested in the average number of steps needed to sort a tabloid of a given shape. Therefore, fix $T \in \op{T}(\lambda)$ and let $r(T,x)$ be the length of the forward slide $\sigma_x$. As $\sigma_x$ can be decomposed into $r(T,x)$ exchanges, the total number of exchanges needed to transform $T$ into a standard Young tableau is
\begin{align*}
r(T) \defeq \sum_{x \in \lambda} r(T,x).
\end{align*}
The average number of needed exchanges is
\begin{align*}
C(\lambda) \defeq \frac{1}{n!} \sum_{T \in \op{T}(\lambda)} r(T).
\end{align*}
We call the number $C(\lambda)$ the complexity of the Novelli--Pak--Stoyanovskii algorithm on the shape $\lambda$. The surprising observation by Krattenthaler and M\"uller can now be stated as $C(\lambda) = C(\lambda')$. In \cite{NeuSul2013} two related objects are considered. The first of these is the drop function which counts the number of tabloids in which the entry $k$, where $1 \leq k \leq n$, is moved to the cell $x \in \lambda$ by its own forward slide. We denote the drop function by
\begin{align*}
d_{\lambda}(k, x) \defeq \#\{ T \in \op{T}(\lambda) : \text{there exists } z \in \lambda \text{ such that } T(z) = k, T_z(x) = k \}.
\end{align*}
The second kind of objects are exchange numbers which count how often two entries are exchanged. That is, given entries $k$ and $l$, where $1 \leq k < l \leq n$, we are interested in the number of tabloids such that the exchange $(k,l)$ appears in the decomposition of a maximal forward slide $\sigma_z$ into exchanges, for some cell $z \in \lambda$. Formally, we define the exchange numbers as
\begin{align*}
\varepsilon_{\lambda}(k,l) \defeq \#\{ T \in \op{T}(\lambda) : k < l \text{ and there exists } z \in \lambda \text{ such that } T(z) = l, \sigma_z(k) \neq k \}.
\end{align*}
For example, in Figure \ref{fig:sort} the entries 5, 9 and 11 all drop to the cell $(2, 4)$. The entry 1 is exchanged with the entries 7 and 11, and the tabloid contributes ten exchanges to the overall complexity $C(\lambda)$, where $\lambda = (4,4,3)$.


\section{Formulae} \label{sec:formulae}

In this section we start out with some preparations we need to define signed exit numbers. Afterwards, we prove formulae for the complexity, the drop function and the exchange numbers, which imply their symmetry in $\lambda$ and $\lambda'$. We shall derive them from Theorem \ref{thrm:exit} which is a result on signed exit numbers.

\bigskip
Let $k$ and $l$ be two entries such that $1 \leq k < l \leq n$ and $x, y \in \lambda$ two cells. We are interested in the number of times $k$ and $l$ are exchanged at the positions $x$ and $y$ during the sorting of all tabloids. Suppose that $T^{-1}(l) = z$, then this happens if and only if $\sigma_z$ moves the entry $k$ from the cell $x$ to the cell $y$. To make this precise, let $z_1 \prec \dots \prec z_n$ be the cells of $\lambda$. We define the local exchange numbers as
\begin{align*}
\varepsilon_{\lambda}(k, l, x, y)
&\defeq \#\{T \in\op{T}(\lambda) : k < l \text{ and } T(z_i) = l, T_{z_{i+1}}(x) = k, T_{z_{i}}(y) = k \text{ for some } i\}.
\end{align*}
It is immediate from their definition that $\varepsilon(k,l,x,y) = 0$ unless $x$ and $y$ are neighbouring cells. More precisely, since $k < l$ the cell $x$ must be a bottom of right neighbour of $y$. An interesting property of local exchange numbers is that they do not depend on the larger argument $l$.

\begin{mylem}\label{lem:locex} Let $k, l, m, n \in \mathbb N$ such that $1 \leq k < l, m \leq n$, let $\lambda \vdash n$ be a partition and $x, y \in \lambda$ be two cells. Then
\begin{align*}
\varepsilon_{\lambda}(k,l,x,y) = \varepsilon_{\lambda}(k,m,x,y).
\end{align*}
\end{mylem}

It is therefore convenient to define $\varepsilon_{\lambda}(k,x,y) \defeq \varepsilon_{\lambda}(k,n,x,y)$ and $\varepsilon_{\lambda}(k)\defeq \varepsilon_{\lambda}(k,n)$. This lemma is proven in \cite[Proposition 4.1]{NeuSul2013}. We shall give a short proof using the Invariance Lemma~\ref{lem:invariance}.

\begin{myproof} Let $T \in \op{T}(\lambda)$, and set $\tilde T \defeq (l,m) \circ T$. Let $z \defeq T^{-1}(l) = \tilde T^{-1}(m)$. By Lemma \ref{lem:invariance} the forward slide $\sigma_z$ moves $k$ from $x$ to $y$ during the sorting of $T$ if and only if $\tilde \sigma_z$ does the same during the sorting of $\tilde T$. Thus, the map $T \mapsto (l,m) \circ T$ provides an involution on $\op{T}(\lambda)$ mapping tabloids which contributing to $\varepsilon_{\lambda}(k,l,x,y)$ to those which contribute to $\varepsilon_{\lambda}(k,m,x,y)$.
\end{myproof}

\bigskip
Following \cite{NeuSul2013}, we define the signed exit numbers as
\begin{align*}
 \Delta_{\lambda}(k,x) \defeq \sum_{y \in N_{\lambda}^{-}(x)} \varepsilon_{\lambda}(k,x,y) - \sum_{y \in N_{\lambda}^{+}(x)} \varepsilon_{\lambda}(k,y,x).
\end{align*}
While this definition is not very intuitive at first glance, the signed exit numbers allow for a unified approach to the exchange numbers and the drop function. For $k \in \mathbb N$ and a cell $x \in \lambda$ we define
\begin{align*}
f_{\lambda}(k, x) \defeq \#\{ T \in \op{SYT}(\lambda) : T(x) = k \}.
\end{align*}

\begin{mythrm}[Signed Exit Numbers]\label{thrm:exit} Let $k, n \in \mathbb N$ such that $1 \leq k < n$, let $\lambda \vdash n$ be a partition and $x \in \lambda$ be a cell. Then we have the recursion
\begin{align} \label{eq:exitrec}
 (n-k)\, \Delta_{\lambda}(k,x) = (n-1)! - \frac{n! \, f_{\lambda}(k, x)}{f_{\lambda}} + \sum_{l = 1}^{k - 1} \Delta_{\lambda}(l, x).
\end{align}
Solving the recursion we obtain
\begin{align} \label{eq:exit}
\Delta_{\lambda}(k, x) &= \frac{n!}{(n - k)(n - k + 1)} \left( 1 - \frac{f_{\lambda}(k, x)}{f_{\lambda}}(n - k + 1) - \sum_{l = 1}^{k - 1} \frac{f_{\lambda}(l, x)}{f_{\lambda}} \right).
\end{align}
\end{mythrm}
\begin{myproof} The recursion \eqref{eq:exitrec} was deduced in \cite[Theorem 5.3]{NeuSul2013} in a slightly more general setting. Its proof makes use of Lemma \ref{lem:locex}.

One can derive \eqref{eq:exit} from \eqref{eq:exitrec} by induction on $k$ where the base of the induction and the induction step are both covered by the same computation:
\begin{align*}
\Delta_{\lambda}(k, x)
&= \frac{(n - 1)!}{n - k} - \frac{n! \, f_{\lambda}(k, x)}{(n - k)\, f_{\lambda}} + \frac{1}{n - k} \sum_{l = 1}^{k - 1} \Delta_{\lambda}(l, x) \\
&= \frac{(n - 1)!}{n - k} - \frac{n! \, f_{\lambda}(k, x)}{(n - k)\, f_{\lambda}} + \frac{n!}{n - k} \sum_{l = 1}^{k - 1} \frac{1}{(n - l)(n - l + 1)} \\
&\quad - \sum_{l = 1}^{k - 1} \frac{n!}{(n - k)(n - l)} \frac{f_{\lambda}(l, x)}{f_{\lambda}} - \sum_{l = 1}^{k - 1} \sum_{i = 1}^{l - 1} \frac{n!}{(n - k)(n - l)(n - l + 1)} \frac{f_{\lambda}(i, x)}{f_{\lambda}} \\
&= \frac{(n - 1)!}{n - k} - \frac{n! \, f_{\lambda}(k, x)}{(n - k)\, f_{\lambda}} + \frac{n!}{n - k} \left( \frac{1}{n - k + 1} - \frac{1}{n} \right) \\
&\quad - \sum_{l = 1}^{k - 1} \frac{n!}{(n - k)(n - l)} \frac{f_{\lambda}(l, x)}{f_{\lambda}} - \sum_{i = 1}^{k - 2} \sum_{l = i + 1}^{k - 1} \frac{n!}{(n - k)(n - l)(n - l + 1)} \frac{f_{\lambda}(i, x)}{f_{\lambda}} \\
&= \frac{n!}{(n - k)(n - k + 1)} - \frac{n! \, f_{\lambda}(k, x)}{(n - k)\, f_{\lambda}} - \sum_{l = 1}^{k - 1} \frac{n!}{(n - k)(n - l)} \frac{f_{\lambda}(l, x)}{f_{\lambda}} \\
&\quad - \sum_{i = 1}^{k - 2} \frac{n!}{n - k} \frac{f_{\lambda}(i, x)}{f_{\lambda}} \left( \frac{1}{n - k + 1} - \frac{1}{n - i} \right) \\
&= \frac{n!}{(n - k)(n - k + 1)} \left( 1 - \frac{f_{\lambda}(k, x)}{f_{\lambda}}(n - k + 1) - \sum_{l = 1}^{k - 1} \frac{f_{\lambda}(l, x)}{f_{\lambda}} \right)
\end{align*}
\end{myproof}


 Note that the recursion \eqref{eq:exitrec} cannot be solved in the setting of \cite{NeuSul2013}. In contrast to the more general approach we obtain the following formula for the drop function.

\begin{mythrm}[Drop Function]\label{thrm:drop} Let $k, n \in \mathbb N$ such that $1 \leq k \leq n$, let $\lambda \vdash n$ be a partition and $x \in \lambda$ be a cell. Then
\begin{align}
d_{\lambda}(k, x)
\label{eq:drop}
&= \frac{n!}{n - k + 1} \left( 1 - \sum_{l = 1}^{k - 1} \frac{f_{\lambda}(l, x)}{f_{\lambda}} \right) \\
\label{eq:drop2}
&= \frac{\prod_{y\in\lambda} h_{\lambda}(y)}{n - k + 1} \sum_{l = k}^{n} f_{\lambda}(l, x).
\end{align}
\end{mythrm}

\begin{myproof} An entry $k$ drops to $x$ if it starts there or is moved there by an exchange with a smaller entry and is not moved away by an exchange with a smaller entry. Thus,
\begin{align} \label{eq:droprec}
d_{\lambda}(k, x)
= (n - 1)! + \sum_{l = 1}^{k - 1} \Delta_{\lambda}(l,x)
= (n-k)\, \Delta_{\lambda}(k,x) + \frac{n!\, f_{\lambda}(k,x)}{f_{\lambda}}.
\end{align}
Substituting \eqref{eq:exit} in \eqref{eq:droprec} yields \eqref{eq:drop}.
Now, \eqref{eq:drop2} follows from \eqref{eq:drop} using $\sum_{k = 1}^{n} f_{\lambda}(k, x) = f_{\lambda}$ and the hook-length formula \eqref{eq:HLF}.
\end{myproof}

\bigskip
Note that \eqref{eq:drop2} can be viewed as a generalised hook-length formula. In the case $n = 1$ it reduces to
\begin {align*}
 n! = n\, (n-1)! = n\, d_{\lambda}(1,x) = \frac{n}{n} \prod_{y \in \lambda} h_{\lambda}(y) \sum_{l = 1}^{n} f_{\lambda}(l,x) = \prod_{y \in \lambda} h_{\lambda}(y) \cdot \#\op{SYT}(\lambda).
\end {align*}


Next, we obtain a description of the exchange numbers.

\begin{mythrm}[Exchange Numbers]\label{thrm:exchange} Let $k, n \in \mathbb N$ such that $1 \leq k < n$, and $\lambda \vdash n$ be a partition. Then we have the recursion
\begin{align}
\label{eq:recexchange}
(n - k)\, \varepsilon_{\lambda}(k) = (n - 1)! \sum_{x \in \lambda} h'(x) + \sum_{l = 1}^{k - 1} \varepsilon_{\lambda}(l) - \frac{n!}{f_{\lambda}} \sum_{x \in \lambda} h'(x)\, f_{\lambda}(k, x).
\end{align}
Furthermore,
\begin{align} \label{eq:exchange}
\varepsilon_{\lambda}(k) &= \frac{n!}{(n - k)(n - k + 1)} \sum_{x \in \lambda} h'(x)
\left( 1 - \frac{f_{\lambda}(k, x)}{f_{\lambda}}(n - k + 1) - \sum_{l = 1}^{k - 1} \frac{f_{\lambda}(l, x)}{f_{\lambda}} \right).
\end{align}
\end{mythrm}

\begin{myproof} The recursion was deduced in \cite[Theorem 4.4]{NeuSul2013} (again in a slightly more general setting) but both claims now follow directly from Theorem \ref{thrm:exit} and the fact that
\begin{align*} 
 \varepsilon_{\lambda}(k)
 = \sum_{x \in \lambda} \sum_{y \in N_{\lambda}^{-}(x)} \varepsilon_{\lambda}(k, x, y)
 = \sum_{x \in \lambda} h'(x) \Delta_{\lambda}(k,x).
\end{align*}
\end{myproof}

\bigskip

For $n \in \mathbb N$ let $H_n \defeq \sum_{k = 1}^n \frac{1}{k}$ denote the $n$-th harmonic number.

\begin{mythrm}[Complexity]\label{thrm:complexity} Let $n \in \mathbb N$ and $\lambda \vdash n$ be a partition. Then we have
\begin{align}
\label{eq:comp1}
C(\lambda)
&= \left( \sum_{i \geq 1} \binom{\lambda_i}{2} + \sum_{i \geq 1} \binom{\lambda'_i}{2}\right) (H_n - 1) - \sum_{x \in \lambda} \sum_{k = 1}^{n - 1} h'(x) \frac{f_{\lambda}(k, x)}{f_{\lambda}} H_{n - k} \\
\label{eq:comp2}
&= \sum_{x \in \lambda} \sum_{k = 1}^{n - 1} h'(x) \frac{f_{\lambda}(k, x)}{f_{\lambda}} (H_{n} - H_{n - k} - 1).
\end{align}
\end{mythrm}

\begin{myproof} On the one hand, the complexity is given in terms of exchange numbers as
\begin{align} \label{eq:sumex}
C(\lambda) = \frac{1}{n!} \sum_{k = 1}^{n - 1} (n - k)\, \varepsilon_{\lambda}(k).
\end{align}
On the other hand, the complexity can be expressed in terms of the drop function as
\begin{align}
\label{eq:sumd1}
C(\lambda)
&= \frac{1}{n!} \sum_{x \in \lambda} \sum_{k = 1}^{n} h'(x) (d_{\lambda}(k, x) - (n - 1)!) \\
\label{eq:sumd2}
&= \frac{1}{n!} \sum_{x \in \lambda} \sum_{k = 1}^{n} h'(x) \left( d_{\lambda}(k, x) - \frac{n!\, f_{\lambda}(k, x)}{f_{\lambda}} \right),
\end{align}
where \eqref{eq:sumd1} counts how often an entry is exchanged with a smaller one (i.e., it drops), and \eqref{eq:sumd2} counts how often an entry is exchanged with a larger one.

There are various possibilities to prove the claims. Substitution of either
\eqref{eq:exchange} in \eqref{eq:sumex} or \eqref{eq:drop} in
\eqref{eq:sumd1} yields \eqref{eq:comp1}. Inserting \eqref{eq:drop2} into
\eqref{eq:sumd2} we obtain \eqref{eq:comp2}.
Moreover, \eqref{eq:comp1} can
easily be transformed into \eqref{eq:comp2} and vice versa using
\begin{align*}
\sum_{x \in \lambda} \sum_{k = 1}^{n} h'(x) \frac{f_{\lambda}(k, x)}{f_{\lambda}}
= \sum_{x \in \lambda} h'(x)
= \sum_{i \geq 1} \binom{\lambda_i}{2} + \sum_{i \geq 1} \binom{\lambda'_i}{2}.
\end{align*}
\end{myproof}


As promised, the symmetry conjecture is now an immediate consequence. For any cell $x = (i,j) \in \lambda$ set $x' \defeq (j,i) \in \lambda'$. For any tabloid $T \in \op{T}(\lambda)$ define a tabloid $T' \in \op{T}(\lambda')$ by letting $T'(x') \defeq T(x)$. For any hook tableau $H \in \op{H}(\lambda)$ define a hook tableau $H' \in \op{H}(\lambda')$ by $H'(x') \defeq -H(x)$. The induced maps $\lambda \to \lambda'$, $\op{T}(\lambda) \to \op{T}(\lambda')$ and $\op{H}(\lambda) \to \op{H}(\lambda')$ are obviously bijective and involutive in the sense that $(x')'=x$, $(T')' = T$, and $(H')' = H$. Moreover, $T \in \op{SYT}(\lambda)$ if and only if $T' \in \op{SYT}(\lambda')$.

\begin{mycor}[Symmetry]\label{cor:symm} Let $k,n \in \mathbb N$ such that $1 \leq k \leq n$, let $\lambda \vdash n$ be a partition and $x \in \lambda$ be a cell. Then
\begin{align*}
C(\lambda) = C(\lambda'),\quad
d_{\lambda}(k,x) = d_{\lambda'}(k,x'), \quad
\varepsilon_{\lambda}(k) = \varepsilon_{\lambda'}(k),\quad \text{ and }\quad
\Delta_{\lambda}(k,x) = \Delta_{\lambda'}(k,x').
\end{align*}
\end{mycor}

\begin{myproof} The assertions follow from Theorems \ref{thrm:exit}, \ref{thrm:drop}, \ref{thrm:exchange} and \ref{thrm:complexity} and the fact that $f_{\lambda}(k,x) = f_{\lambda'}(k,x')$ by virtue of the bijection $T \mapsto T'$.
\end{myproof}

However, note that $\varepsilon_{\lambda}(k,x,y) = \varepsilon_{\lambda'}(k, x', y')$ is false in general.


\bigskip
In the rest of this section we explain how the complexity $C(\lambda)$ of the Novelli--Pak--Stoyanovskii algorithm on a given shape $\lambda \vdash n$ can actually be computed.
 
The most primitive approach is certainly to sort all $n!$ tabloids, counting the exchanges. Theorem~\ref{thrm:complexity} allows for a somewhat better method. That is, now we ``only'' need to consider all standard Young tableaux of shape $\lambda$, and count how often each entry occurs in each cell, in order to get a hold on the numbers $f_{\lambda}(k, x)$. In some special cases these numbers can be given a closed form. In general the cost of computing $f_{\lambda}(k,x)$ ``only'' depends on the number of sub-partitions of $\lambda$. To see this, decompose a tabloid $T$ with $T(x) = k$ into the parts $U$ with entries smaller than $k$ and $V$ with entries larger than $k$. Then $U$ is a standard Young tableau of shape $\mu$ for some $\mu \subseteq \lambda$ while $V$ corresponds to a skew standard Young tableaux. Hence, the numbers $f_{\lambda}(k,x)$ can be computed by counting standard Young tableaux and skew tableaux whose shape depends on a suitable $\mu$. This can be done using the hook-length formula \eqref{eq:HLF} for standard Young tableaux and Aitken's determinant formula \cite{Aitken1943} for the skew tableaux.

\section{Bijections} \label{sec:bijections}

We are now going to give bijective proofs of the formulae and symmetry results in the last section.

\bigskip
Let $\op{T}(\lambda, k \to x)$ denote the set of all tabloids of shape $\lambda$ such that the entry $k$ drops to the cell $x$ during the application of the sorting algorithm. Moreover, let $\op{SYT}(\lambda,x \geq k)$ be the set of standard Young tableaux of shape $\lambda$ such that the entry of the cell $x$ is at least $k$. Theorem \ref{thrm:drop} \eqref{eq:drop2} suggests we should look for a bijection
\begin{align*}
 \Psi: \op{T}(\lambda, k \to x) \times \{k, \dots, n\} \to \op{H}(\lambda) \times \op{SYT}(\lambda, x \geq k).
\end{align*}
Indeed, given a tabloid $T \in \op{T}(\lambda, k \to x)$ we notice that the bijection \eqref{eq:nps} maps $T$ to a pair $\Phi(T) = (H, U)$ of a hook tableau $H$ and a standard Young tableau $U \in \op{SYT}(\lambda, x \geq k)$. This is true because $k$ drops to the cell $x$, and can only be moved away by a larger entry afterwards. Now, consider a pair $(T, l)$ of a tabloid $T \in \op{T}(\lambda, k \to x)$ and an integer $k \leq l \leq n$. If we apply $\Phi$ to the tabloid $(kl) \circ T$ then by the Invariance Lemma \ref{lem:invariance} we will again obtain a hook tableau $H$ and a standard Young tableau $U \in \op{SYT}(\lambda, x \geq k)$.

\begin{mythrm}\label{thrm:dropbij} Let $k, n \in \mathbb N$ such that $1 \leq k \leq n$, $\lambda \vdash n$ be a partition and $x \in \lambda$ be a cell. Then the map 
\begin{align*}
 \Psi: \op{T}(\lambda, k \to x) \times \{k, \dots, n\} &\to \op{H}(\lambda) \times \op{SYT}(\lambda, x \geq k), \\
 (T,l) &\mapsto \Phi((kl) \circ T),
\end{align*}
where $\Phi$ is given by the Novelli--Pak--Stoyanovskii algorithm, is a bijection.
\end{mythrm}

\begin{myproof} The remarks above show that $\Psi$ is well-defined. We prove the claim by constructing an inverse map. Therefore, let $(H,U) \in \op{H}(\lambda) \times \op{SYT}(\lambda, x \geq k)$. Since the map $\Phi: \op{T}(\lambda) \to \op{H}(\lambda) \times \op{SYT}(\lambda)$ is a bijection, there is a well-defined tabloid $T \defeq \Phi^{-1}(H,U)$. If $T \in \op{T}(\lambda, k\to x)$ then set $\Psi^{-1}(H, U) \defeq (T, k)$. In this case we have $\Psi^{-1} \circ \Psi (T, k) = (T, k)$ and $\Psi \circ \Psi^{-1}(H,U) = (H, U)$.

If $T \notin \op{T}(\lambda, k \to x)$ then we have to determine an integer $l$, where $k < l \leq n$, such that $k$ drops to $x$ during the sorting of $(kl) \circ T$. To this end, let $I(x) \subseteq \lambda$ be the set of cells which lie to the right and below $x$. That is, $z \in I(x)$ if and only if $z \in \lambda$ and $z \geq x$ with respect to the component-wise order.

We consider the tabloid $T_x$. If $T_x(z) \geq k$ for all cells $z \in I(x)$ then we set $l \defeq T(x)$, and claim that $(kl) \circ T \in \op{T}(\lambda, k \to x)$. To see this, let $\tilde T \defeq (kl) \circ T$. By the Invariance Lemma \ref{lem:invariance} we have $\tilde T_x(z) \geq k$ for all $z \in I(x)$. Since $\tilde T_x$ is sorted on $I(x)$ it follows that $\tilde T_x(z) > k$ for all $z \in I(x) \setminus \{x\}$ and $\tilde T_x(x) = k$. But this means that the entry $k$ has not moved at all but dropped to its initial cell $x$.

Otherwise, there is a $z \in I(x)$ such that $T_x(z) < k$. Indeed, we must have $T_x(x) < k$ since $T_x$ is sorted on $I(x)$. In this case let $y$ be the maximal cell with respect to the order $\prec$ used in the Novelli--Pak--Stoyanovskii algorithm such that $T_y(z) \geq k$ for all $z \in I(x)$. Such a $y$ must exist because $U(x) \geq k$. We set $l \defeq T(y)$ and claim that $\tilde T \defeq (kl) \circ T \in \op{T}(\lambda, k \to x)$ as before. By the Invariance Lemma \ref{lem:invariance}, $y$ is the maximal cell with respect to $\prec$ such that $\tilde T_y(z) \geq k$ for all $z \in I(x)$. This means that the entry $\tilde T(y)$ must drop to a cell in $I(x)$. But $\tilde T(y) = k$ is the smallest entry of $\tilde T_y$ among the cells of $I(x)$. Since $\tilde T_y$ is ordered on $I(x)$, the entry $k$ must drop to the cell $x$.

Once we have found a suitable $l$ as described above, we set $\Psi^{-1}(H,U) \defeq ((kl) \circ T, l)$. Clearly, we have $\Psi \circ \Psi^{-1}(H,U) = (H,U)$.

Conversely, let $T \in \op{T}(\lambda, k \to x)$ and $k < l \leq n$. Set $\tilde T \defeq (kl) \circ T$ and $y \defeq T^{-1}(k) = \tilde T^{-1}(l)$. Then $y = x$ or $y$ is the maximal cell of $\lambda$ with respect to $\prec$ such that $\tilde T_y(z) \geq k$ for all $z \in I(x)$. This once again follows from the Invariance Lemma \ref{lem:invariance}.

Hence, $\Psi^{-1} \circ \Psi(T,l) = (T,l)$, and the proof is complete.
\end{myproof}

\begin{figure}
 \begin{center}
  \includegraphics[height=30mm]{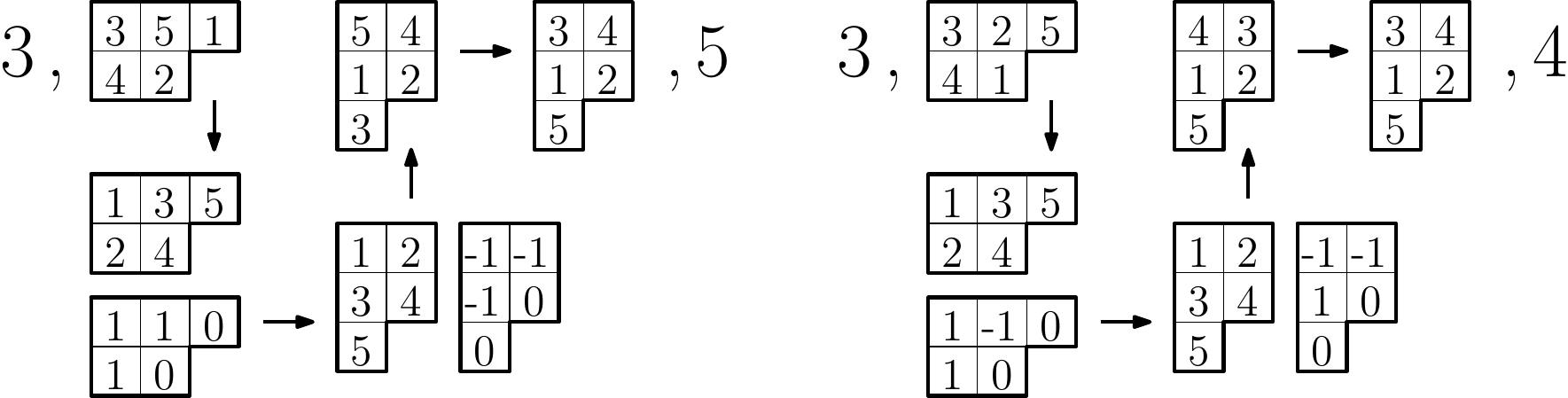}
  \caption{Two examples where $\lambda = (3,2)$, $k=3$ and $x=(1,2)$.}
  \label{fig:dropbij}
 \end{center}
\end{figure}

\bigskip
Note that $\Psi$ induces a bijection
\begin{align}
\op{T}(\lambda, k \to x) \times \{k, \dots, n\} &\to \op{T}(\lambda', k \to x') \times \{k, \dots, n\}, \notag \\ \label{eq:dropbij}
(T,l) \overset{\Psi_{\lambda}}{\longmapsto} (H, U) &\mapsto (H', U') \overset{\Psi_{\lambda'}^{-1}}{\longmapsto} (V, m),
\end{align}
which is an involution in the obvious sense. Unfortunately, it does not seem to be easy to eliminate the auxiliary integers. As Figure \ref{fig:dropbij} shows, even in simple cases distinct tabloids may be mapped to the same tabloid with different integer labels.

Using the ideas from Theorem \ref{thrm:dropbij}, we want to find a bijective proof of the formula \eqref{eq:exchange} for the number of exchanges of $k$ with a larger entry. Suppose during the sorting of the tabloid $T \in \op{T}(\lambda)$ the entries $k$ and $l$ with $k < l$ are exchanged. We identify this exchange with the pair $(T,l)$. Let $\op{Ex}(\lambda, k)$ denote the set of all such exchanges, i.e. $\# \op{Ex}(\lambda, k) = (n - k)\, \varepsilon_{\lambda}(k)$.

Now, fix a tabloid $U$. Suppose that $k$ drops to $x$ during the sorting of $U$, and is afterwards moved to the final cell $y$. Clearly,
\begin{align*}
 \# \{ (T,l) \in \op{Ex}(\lambda, k) : T = U \} = h'(x) - h'(y).
\end{align*}
We denote this cardinality by $e(U,k)$. Thus, we obtain a bijection
\begin{align*}
 \op{Ex}(\lambda, k) \to \{ (T,i) : T \in \op{T}(\lambda), 1\leq i \leq e(T,k) \}, 
\end{align*}
for example by ordering the exchanges $(T, l_i)$, $1 \leq i \leq e(T,k)$ with respect to the order in which they occur during the sorting. Thereby, letting 
\begin{align*}
 A(\lambda, k) \defeq \bigcup_{x \in \lambda} \op{T}(\lambda, k \to x) \times \{ 1, \dots, h'(x) \},
\end{align*}
we have an inclusion $\iota: \op{Ex}(\lambda, k) \to A(\lambda, k)$. Moreover, set
\begin{align*}
 B(\lambda, k) \defeq \op{H}(\lambda) \times \bigcup_{y \in \lambda} \{T \in \op{SYT}(\lambda) : T(y) = k \} \times \{1, \dots, h'(y)\}.
\end{align*}

\begin{myprop}\label{prop:exchangebij} Let $k,n \in \mathbb N$ such that $1 \leq k \leq n$, and $\lambda \vdash n$ be a partition. Then the map
\begin{align*}
\psi: A(\lambda, k) \setminus \op{Ex}(\lambda, k) \to B(\lambda, k)
\end{align*}
defined by $(T, i) \mapsto (H, U, i-e(T,k))$, where $(H,U) \defeq \Phi(T)$ is given by the Novelli--Pak--Stoyanovskii algorithm, is a bijection.
\end{myprop}

\begin{myproof} Let $(T,i)$ be a pair of a tabloid and an integer, and let $x$ be the cell to which $k$ drops during the sorting of $T$. Set $(H,U) \defeq \Phi(T)$ and $y \defeq U^{-1}(k)$. Then the inequalities
\begin{align*}
 e(T, k) = h'(x) - h'(y) < i \leq  h'(x) \quad \text{and} \quad
 0 < i - e(T,k) \leq h'(y)
\end{align*}
are equivalent. Thus, $(T,i)$ lies in $A(\lambda, k)$ but not in $\op{Ex}(\lambda, k)$ if and only if $(H,U,i-e(T,k))$ lies in $B(\lambda, k)$. It follows that $\Psi$ is well defined and a bijection.
\end{myproof}


Considering only cardinalities, where we use Theorem \ref{thrm:dropbij} to specify the cardinality of $A(\lambda, k)$, Proposition \ref{prop:exchangebij} reduces to \eqref{eq:exchange}. Lastly, we want to give a bijective proof of $C(\lambda) = C(\lambda')$. We do this by indicating a bijection
\begin{align*}
 \tilde \psi : \op{Ex}(\lambda, k) \times \{ k, \dots, n \} &\to \op{Ex}(\lambda', k) \times \{k, \dots, n\}.
\end{align*}
The bijections $\Psi$ and $\psi$ of Theorem \ref{thrm:dropbij} and Proposition \ref{prop:exchangebij} provide a bijection
\begin {align*}
 f: (\op{Ex}(\lambda, k) \cup B(\lambda, k)) \times \{k \dots, n\} \to (\op{Ex}(\lambda', k) \cup B(\lambda', k)) \times \{k \dots, n\}.
\end {align*}
However, we have an additional bijection $g: B(\lambda, k) \to B(\lambda', k)$ given by $(H,T,i) \mapsto (H',T',i)$ at our disposal. Consider a sequence
\begin{align*}
 (a,l) \overset{f}{\longmapsto} (b_1',l_1') \overset{g}{\longmapsto} (b_1,l_1) \overset{f}{\longmapsto} (b_2',l_2') \overset{g}{\longmapsto} (b_2,l_2) \overset{f}{\longmapsto} \dots
\end{align*}
where $a \in \op{Ex}(\lambda, k)$, $b_j' \in B(\lambda', k)$ and $b_j \in B(\lambda, k)$. Since $f$ and $g$ are bijections, no pair $(b_j',l_j')$ or $(b_j ,l_j)$ can recur. Since the considered sets are finite the sequence must end at some point with a pair $(a',l')$ where $a' \in \op{Ex}(\lambda',k)$. Hence, by setting $\tilde \psi(a,l) \defeq (a',l')$ we have found the desired bijection.

\section*{Acknowledgements}\label{sec:ack}

The author wants to thank Christian Krattenthaler and Cristoph Neumann for their helpful comments.

\bibliographystyle{alpha}
\bibliography{robin}
\label{sec:biblio}

\end{document}